\documentclass{article}
\usepackage{amsmath,amssymb}
 
 \newcommand{\ZZ}{\mathbb{Z}}
 \newcommand{\NN}{\mathbb{N}}
\newcommand{\RR}{\mathbb{}}

 \def\auts{cellular automata}
\newcommand{\aut}{{cellular automaton}}
 \newcommand{\autb}{{cellular automaton}} 
\newcommand{\autsb}{{cellular automata }}
 \newcommand{\az}{{A^\ZZ}}

\newtheorem{thm}{Theorem}[section]
 \newtheorem{pro}{Proposition}[section]

 \newtheorem{lem}{Lemma}[section]
\newtheorem{cor}{Corollary}[section]

\begin{document}
\title { Some properties of cellular automata with equicontinuity points.}
\author{ F. BLANCHARD and P. TISSEUR \thanks{ Institut de Math\'ematiques de
Luminy - UPR 9016 du CNRS}} 
\sloppy 
\date{\ }
\maketitle  
\begin{abstract}
We investigate topological and ergodic properties of cellular automata
having equicontinuity points. In this class surjectivity on a transitive SFT
implies existence of a dense set of periodic points. Our main result is that
under the action of such an automaton any shift--ergodic measure converges in
Ces\`aro mean, assuming equicontinuity points have measure 1; the limit measure
is described by a formula and some of the properties of its topological support
are given.
\end{abstract}

\section{Introduction.}

Compared with their topological dynamics, the ergodic theory of cellular
automata is still in its infancy. One of the main reasons is that few invariant
measures are known, if any, for any given cellular
automaton. In this article for a family of CA  defined by a
topological property we give a rather simple construction of measures that are
invariant both for the shift and the automaton; we also show that when an
automaton belonging to this family is onto, there is a dense set of periodic
points for its action.

The property of having equicontinuity points was
first considered for \autsb by Gilman \cite{Gil87} in relation with Wolfram's
empirical classification \cite{Wo86}. Afterwards this property was used by
K\accent23urka  \cite{Ku94} as a basic element of his topological
classification of CA according to their local behaviour. He distinguishes
four classes: equicontinuous automata (E1); those that have equicontinuity
points without being equicontinuous (E2); and two other classes of CA, all of
them sensitive to initial conditions.

In \cite{Li84} Lind finds the exact Ces\`aro mean limit of the images of
Bernoulli measures by a simple additive \aut. His result was generalized
in \cite{MM} and later extended to Markov measures
under the action of a larger class of additive \auts, all of them without
equicontinuity points, by Ferrari,
Maass, Martinez and Ney \cite{FAM}. Boyle and Kitchens \cite{BK} also
proved that periodic
points are dense for left-- or right--closing automata; this class has only a
small overlap with E2.

These results concern CA acting on the full set of
configurations $A^\ZZ$. Our setting is slightly larger: we consider a \autb
$F$ acting on a subshift $X \subset A^\ZZ$, that is,
$F(X) \subset X$, and the equicontinuity points we consider are those of the
dynamical system $(X,F)$.

The article is devised as follows. After the Introduction, Section 2 is
devoted to general  definitions, in particular that of blocking words, which is
essential for CA  having equicontinuity points. Section 3 contains the results.
The first one is topological but we prove it with the help of
Poincar\'e's recurrence theorem: assuming that the automaton $F$ acts
surjectively on a transitive subshift of finite type $X$ and
has equicontinuity points, $X$ contains a dense set of
$F$--periodic points. Our main result is purely ergodic:
if a blocking word $(X,F)$ has positive measure for a shift--ergodic measure
$\mu$, then the images of $\mu$ under the powers of $F$
converge in Ces\`aro mean. The limit $\mu_c$ is of course $F$-- and
$\sigma$--invariant; it is given by a formula, which nevertheless leaves open
several interesting questions. Then we examine properties of the topological
support of $\mu_c$, and give several examples  of \autsb to which our results
can be applied; we finish with some open questions.

\section{Definitions and background.}

\subsection{Dynamical systems, measures, cellular automata.} $\mbox{ }$

A topological dynamical system $(X,T)$ consists of a compact metric space $X$
and a continuous self--map $T$.

A point $x\in X$ is said to be an equicontinuity point, or to be Lyapunov
stable,
if for any $\epsilon >0$, there exists $\eta >0$ such that if $d(x,y)\le
\eta$ one has $d(T^i(x),T^i(y))\le \epsilon$ for any integer $i>0$. When all
points of $(X,T)$ are equicontinuity points $(X,T)$ is said to be
equicontinuous: since $X$ is compact an equicontinuous system is
uniformly equicontinuous.

The dynamical systems in this article are all defined on a symbolic space. Let
$A$ be a finite alphabet. Define $A^*$ to be the set of all finite
concatenations of letters of $A$, called words; the length of the word $u$ is
denoted by $\vert u\vert $. A language $L$ is a subset of $A^*$.

The set $A^\ZZ$ of bi--infinite sequences on $A$ is endowed with the product
topology, associated with the distance $d(x,y)=2^{-i}$ where
$i=\min\{\vert j\vert x(j)\ne y(j)\}$; the shift $\sigma$: $\sigma (x)=
(x_{i+1})_{i\in\ZZ}$ is a homeomorphism on $A^\ZZ$.
Given a word $u$ and an integer $t$, the clopen set $[u]_t=\{x\in
\az : x_t=u_1\ldots ;x_{t+\vert u\vert -1}=u_{\vert u\vert}\}$ is called a
cylinder set. When $x\in A^\ZZ$ and $p \le q$ are two integers, put
$x(p,q)=x_p\ldots x_q$. A sequence $(S_i)_{i\in\NN}$ is said to be
ultimately periodic if there exist two natural integers $p'$ and $p$ such that
$S_{p'+kp+i}=S_{p'+i}$ for $k$, $i\in \NN$.

The dynamical system $(A^\ZZ ,\sigma )$ is called the full shift. A
subshift $X$ is a closed shift--invariant subset of $A^\ZZ$. A transitive
subshift $S$ is one such that for $u, v \in L(S)$ there is $w\in L(S)$ such that
$uwv\in L(S)$; it is strongly mixing if for any $n$ larger than some
$n_0(u,v)$ one can find a word $w$ of length $n$ with the same property. To
any subshift $X$ there corresponds a unique language $L(X)$: it is the set of
all words that are found as  blocks of coordinates of a point
of $X$. Given any subshift $X$ the language $L(X)$ has two general properties:
for $u\in L(X)$, any sequence of consecutive letters of $u$ is also in
$L(X)$; and for any word $v$ in $L(X)$ there are letters $a$ and $b$ in $A$ such
that $avb\in L(X)$). A subshift of finite type $X$ is defined by forbidding a
finite family of words $E$: then $L(X)$ is the smallest language having the two
properties above and such that no word  $u \in L(X)$ is of the form $u=vew$ with
$e\in E$. Transitive subshifts of finite  type have a  dense set of periodic
points.

When a probability measure  $\mu$ on $A^\ZZ$ is shift--invariant, its
topological support $S(\mu )$ is closed invariant, hence a
subshift. On every transitive subshift of finite type one defines a
particular measure $\lambda$ with support $X$ called the Parry measure;
the Parry measure of the full shift is the
uniform measure.

A sequence
$(\mu_n )_{n\in\NN}$ of probability measures on a compact set $K$ is said to
converge vaguely to a limit  $\mu$ if the sequence $\int_{K}fd\mu_n$
tends to $\int_{K}fd\mu $ for any continuous function $f\colon K\to \RR$.
On a subshift $X$ a sequence
$(\mu_n)_{n\in\NN}$ of shift--invariant measures converges vaguely if and only
if for any word $u\in L(X)$ the sequence $\left (\mu_n
([u]_{0})\right )_{n\in\NN}$ converges.


In this article we call \autb (CA for short) a continuous map
$F : X\to X$ defined on a subshift $X\subseteq A^\ZZ$ and commuting with the
shift $\sigma$; we also call \autb the dynamical system $(X,F)$. The
Curtis--Hedlund--Lyndon theorem \cite{DGS75} states that for every \autb $(X,F)$
there is an integer $r$, called the radius of $F$, and a block map $f :
A^{2r+1}\cap L(X)\to A$ such that one has $$ F(x)_i=f(x_{i-r},\ldots ,x_i
,\ldots ,x_{i+r}). $$

If $X$ is  a transitive subshift of finite type, the automaton $F$ acts
surjectively on $X$ if and only if the Parry measure
$\lambda$ is $F$--invariant \cite{CP}.

The set
$W(X,F)=\lim_{n\to\infty}\bigcap_{i=0}^n F^i (X)$ is called the limit set
of the \autb $(X,F)$; of course when $F$ is surjective $W(X,F) = X$.

\subsection{Blocking words and equicontinuous points.}

\begin{ttdefi}
Let $F$ be a \autb with radius $r$ acting on the subshift $X$. A word
$B\in A^{2k+1}$ is called a blocking word for $(X,F)$ if there is an infinite
sequence of words $v_n$, $\vert v_n\vert =2i+1\ge r$ such that for any $x \in
A^\ZZ$ with $x(-k,k)=B$ one has $F^n(x)(-i,i)=v_n$ for $n\in\ZZ^*$.
\end{ttdefi}

In the definition above we do not assume $X$ to be $\sigma$--transitive, but
this condition appears necessary for most of our proofs. Remark that if $B$ is
a blocking word and  $x(-k,k)=B$, then $F^n(x)(-\infty,-i)$ does not depend on
$x(k,+\infty)$, and reversely since $2i+1\ge r$. An occurrence of a
blocking word
in a configuration $x$ completely disconnects coordinates to its right  and left
for the action of the automaton; hence the name.


The two following results are essentially due to K\accent23urka
\cite{Ku94}.

\begin{pro}\label {p5}
Any equicontinuity point of a \autb $(X,F)$ has an occurrence of a blocking
word. Conversely if there exist blocking words, any point with infinitely many
occurrences of a blocking word to the left and right is an equicontinuity
point; if moreover $X$ is transitive for  $\sigma$ equicontinuity
points are dense in $X$. \end{pro}

\begin{ppreuve}
Let $x$ be an equicontinuity
point of $(X,F)$; applying the definition of equicontinuity
points to \autsb, there is an integer
$k$ such that if  $d(x,y)<2^{-k}$, for any  $n$ one has $F^n x(-r,r)=F^n
y(-r,r)$, so that $B= x(-k,k)$ is a blocking word. Conversely let $V(X)$
be the set of all points with infinitely many occurrences of blocking
words to the left and right; when $X$ is transitive $V(X)$ is non--empty, even
dense. Let $x \in V(X)$. To any given $\varepsilon>0$ one associates an integer
$t$ such that $2^{-t}<\varepsilon$. There exist a real number $\eta$ and
integers $t<k$ such that $2^{-k}<\eta$ and the words $x(-k,-t)$ and $x(t,k)$
contain an occurrence of $B$ each. For every point $y$ belonging to
the cylinder set $[x(-k,k)]_{-k}$ one has $F^i(x)(-t,t)=F^i(y)(-t,t)$; since
$\varepsilon$ is chosen arbitrarily one concludes that $x$ is an equicontinuity
point. \end{ppreuve}

For an equicontinuous \aut, there is a natural integer $k$ such that all words
of $L(X)$ with length $2k+1$ are blocking words; thus

\begin{pro}\label{equi1}
The \autb $(X,F)$  is equicontinuous if and only if there are two integers $p$
and $p'$ such that for $x\in X$ the sequence $(F^n(x))_{n\in\NN}$ is ultimately
periodic with period $p$ and preperiod $p'$.
 \end{pro}

\section{Results.}

\subsection{Dense periodic points.}

\begin{pro} \label{dense}
Let $X$ be a transitive subshift of finite type and suppose that the \autb
$(X,F)$ has an equicontinuity point. Then $F$ is onto if and only if it
possesses a dense set of periodic points.  \end{pro}
\begin{ppreuve}  Let $F$ act surjectively on the transitive subshift of finite
type $X$ and suppose it has equicontinuity points. For any word
$v \in L(X)$ we construct a $\sigma$--periodic point  $\bar{u}\in
X$ such that $\bar{u}(k,|v|-1+k) = v$ for some integer $k$, which is also
$F$--periodic; this establishes the density of
$F$--periodic points in $X$.

Fix $v \in L(X)$. By
Proposition \ref{p5} $F$ has a blocking word $B\in L(X)$; as $X$ is
transitive and has a dense set of $\sigma$--periodic points, there is a word
$u=Bwvw' \in L(X)$ such that $\bar{u} \in X$, where $\bar{u}$ is the
periodic point constructed on $u$ and such that an occurrence of $u$ starts at
$0$. The cylinder set $C = [uB]_0$ contains $\bar{u}$, and
$\lambda(C)>0$ if $\lambda$ is the Parry measure of $X$. Since $\lambda$ is
$F$--invariant we apply the Poincar\'e recurrence theorem: there is
$m>0$ such that $\lambda(C\cap F^{-m}C) > 0$; in particular there are a
point $x\in X$ and $q = (\vert uB\vert -1)$ such that $x(0,q) = F^m(x)(0,q) =
uB$.

But $B$ is a blocking word. All the coordinates of $\bar{u}$ coincide with
those of $x$ on the segment $[0,q]$, and since there is an occurrence of $B$ at
the beginning of this segment and one at the end, for any $n>0$ one has $F^nx(i,
q-i) =F^n\bar{u}(i, q-i)$, where $i<\frac{1}{2}\vert B\vert$ is as in the
definition of blocking words. We have thus shown
that $F^nx$ and $F^n\bar{u}$ coincide on a segment of length $q-2i\ge q-|B|=
|u|$, which is greater than or equal to the common $\sigma$--period of $\bar{u}$
and $F^m\bar{u}$: therefore $F^m\bar{u}=\bar{u}$.

The converse is straightforward.
\end{ppreuve}

We have proved this topological result ergodically. There should be a purely
combinatorial proof. The following simple
consequence is known but seems to be nowhere in written form.

\begin{cor}\label{esurj}
A \autb $(X,F)$ is equicontinous and surjective if and only if there is $p>0$
such that any $x\in X$ is periodic of period $p$. \end{cor}
\begin{ppreuve}

By Proposition \ref{equi1}, $F$ being equicontinuous, there is an integer $p'$
such that for any $x\in A^\ZZ$ the sequence $(F^{p'+n}(x))_{n\in\NN}$ is
periodic with period $p$; then any periodic point has period
$p$. By Proposition \ref{dense} the set of periodic points is
dense; as in the proof of this proposition one identifies the block of
coordinates $F^n(x)(-k,k)$ with the corresponding block of a periodic point
with period $p$ for every $n$, and one reaches the conclusion by letting $k$ go
to infinity. 
\end{ppreuve}

\subsection{Ces\`aro mean convergence of measures.}

We start with an easy result on equicontinuous CA. If $\mu$ is a measure on
$A^\ZZ$ and $M$ a Borel set, denote by $$ \mu_n
(M)=\frac{1}{n}\sum_{i=0}^{n-1}\mu \left (F^{-i} (M)\right ) $$ its
Ces\`aro mean of order $n$ with respect to $F$.

\begin{pro}\label{p4}
Let $(X,F)$ be an equicontinous \autb with period $p$ and preperiod $p'$, and
let  $\mu$ be a shift--ergodic measure with support $X$.
Then $\mu$ converges vaguely in Ces\`aro mean to the measure
$\displaystyle{\mu_c=\frac {1}{p}\sum_{i=0}^{p-1}\mu\circ F^{-(i+p')}}$.
\end{pro}

\begin{ppreuve}

It suffices to show that for $u \in L(X)$ the sequence $\left
(\mu_n ([u]_{0})\right )_{n\in\NN}$ converges to the right limit.
By Proposition \ref{equi1} there are $p$ and $p'$ that for any point
$x$, any pair of integers $n$ and $i$ one has
$F^{p'+i+np}(x)=F^{p'+i}(x)$. Thus if $u \in L(X)$ and
$n>p'$ one has

$$
\mu_n ([u]_{0})=\frac{1}{n}\sum_{i=0}^{p'-1}\mu \left ( F^{-i} \left
([u]_{0}\right )\right ) +\frac{1}{n}\sum_{i=p'}^{n-1}\mu \left ( F^{-i}
\left ([u]_{0}\right )\right ). 
$$ 
The first term tends to $0$; using periodicity one gets

$$
\mu_n ([u]_{0})=\frac{1}{p}\sum_{i=0}^{p-1}\mu \left ( F^{-(i+p')} ([u]_{-k}
)\right ). 
$$ 
\end{ppreuve}

\begin{ttdefi}
Let $F$ be a \autb acting on the subshift $X$. A probability measure $\mu$
on $X$
is said to be equicontinuous for $(X,F)$ if the set of equicontinuity points of
$(X,F)$ has measure $1$. 
\end{ttdefi}

\begin{lem}\label {p6}
Let $(X,F)$ be a \autb and $\mu$ be a measure on $X$, ergodic for $\sigma$.
Then the two following properties are equivalent: \itemize
\item{(1)} there exists a blocking word $B$ such that $\mu([B]_0)>0$;
\item{(2)}  $\mu$ is equicontinuous for $(X,F)$. 
\end{lem}

\begin{ppreuve}
(1) $\Rightarrow$ (2): since $\mu$ is $\sigma$--ergodic and $\mu([B]_0)>0$,
almost every point contains infinitely many occurrences of $B$ to the left
and right, so it is an equicontinuity
point by Proposition \ref{p5}. \\ (2) $\Rightarrow$ (1):
again by Proposition \ref{p5}, every equicontinuity
point contains an
occurrence of a blocking word; since the family of blocking words is at
most countable, there is a blocking word $B$ such that $\mu ([B]_{0})>0$.
\end{ppreuve}

\begin{ttdefi}
Given a word $B$, which we shall always suppose to be a blocking word for
$(X,F)$, let $R(k,m,B)$ be the set of all  points of $X$ having at least one
occurrence of $B$ between the coordinates $-m-k$ and $-k$, and another one
between the coordinates $k$ and $m+k$. Whenever there is no ambiguity on $B$ we
denote it by $R(k,m)$. 
\end{ttdefi}

\begin{thm}\label {p7}
Let $(X,F)$ be a \autb and  $\mu$ be a shift--ergodic, equicontinuous measure on
$X$. Then  $\mu$ converges vaguely in Ces\`aro mean under $F$. The limit
$\mu_c$ is $F$-- and $\sigma$--invariant, and for every word $u\in L(X)$ one
has $$ \mu_c([u]_0)=\lim_{m\to\infty} \frac{1}{p(k,m)}\sum_{i=0}^{p(k,m)-1} \mu
\left (R(k,m)\cap F^{-(i+\bar {p}(k,m))} ([u]_{0})\right ). $$ 
\end{thm}

\begin{ppreuve}
It is sufficient to show that for any word $w \in L(X)$, $\vert w\vert =2k+1$,
the sequence $\left (\mu_n ([w]_{-k})\right )_{n\in\NN}$ converges.

By Lemma \ref {p6} there is a blocking word $B$
for $(X,F)$ with $\mu([B]_0)>0$.
The limit of the increasing sequence of sets $\left ( R(k,m)\right
)_{m\in\NN}$ is the set of all points having at least two occurrences of
$B$, one to the left of $-k$ and the other to the right of $k$. Since $\mu$ is
$\sigma$--ergodic the set $V(B)$ of  points having infinitely many occurrences
of $B$ to the right and left has measure $1$. Thus
$\lim_{m\to\infty}\mu ( R(k,m))=1$ and for any integer $k$, any word $w\in
L(X)\cap A^{2k+1}$ one has

$$
\mu_n([w]_{-k})=\lim_{m\to\infty}\frac{1}{n}\sum_{i=0}^{n-1}\mu ( F^{-i} (
[w]_{-k} ) \cap R(k,m)). 
$$

We prove that
$\mu_n ([w]_{-k})$ by using the twofold convergence of the  double sequence
$(\frac{1}{n}\sum_{i=0}^{n-1}\mu ( F^{-i} ( [u]_{-k} ) \cap
R(k,m)))_{m,n\in\NN}$. Indeed since the interval $[0,1]$ in which $\mu$ takes
its values is compact, if $(\frac{1}{n}\sum_{i=0}^{n-1}\mu (F^{-i} ( [u]_{-k} )
\cap R(k,m)))_{m,n\in\NN}$ converges simply as $n\to\infty$ and
uniformly in $n$ as $m\to\infty$, the two limits commute and one obtains
the desired convergence: 
$$
 \lim_{n\to\infty}\lim_{m\to\infty}
\frac{1}{n}\sum_{i=0}^{n-1}\mu ( F^{-i} ( [u]_{-k} ) \cap
R(k,m))=\lim_{n\to\infty}\mu_n([w]_{-k}). 
$$

\medskip

Let us show first that the sequence converges as $n\to \infty$ for fixed $m$.

Let $x$ and $y$ belong to $R(k,m)$: by the definition of $R(k,m)$
there are blocking words to the left of their $-k^{th}$ coordinate and to the
right of their $k^{th}$ coordinate, so that if $y$ is such that
$u=y(-m-k,m+k)=x(-m-k,m+k)$, then for any integer $i$ one has
$F^i(x)(-k,k)=F^i(y)(-k,k)$. In particular if $\bar{u}$ is the periodic point
with period $2m+2k+1$ such that $\bar{u}(-m-k,m+k)=u$, the sequence
$(F^n(\bar{u}))_{n\in\NN} = (F^n(x)(-k,k))_{n\in\NN}$ is ultimately periodic.
Let $p(x,k,m)$ be its period and $p'(x,k,m)$ be its preperiod. Denote by
$p(k,m)$
the least common multiple of the values of $p(x,k,m)$ for $x \in R(k,m)$ and by
$p'(k,m)$ the corresponding integer for $p'(x,k,m)$.

Let $w$ be a word of length $2(k+m)+1$ such that $[w]_{-k-m}\subset R(k,m)$.
For any $x\in [w]_{-k-m}$ and $i$, $j\in \NN$ one has
$F^{p'(k,m)+j+ip(k,m)}(x)(-k,k)=F^{p'(k,m)+j}(x)(-k,k)$; thus for any word $u$
of length $2k+1$ one has 
$$
R(k,m)\cap
F^{-(ip(k,m)+j+p'(k,m)}([u]_{-k})=R(k,m)\cap F^{-(p'(k,m)+j)}([u]_{-k}).
$$
 An argument similar to that of the proof of Proposition \ref{p4} shows that

$$
\hskip -3 true cm \lim_{n\to\infty}\frac{1}{n} \sum_{i=0}^{n-1}\mu \left (F^{-i}
([u]_{-k})\cap R(k,m)\right ) 
$$
$$
\hskip 2 true cm =\frac{1}{p(k,m)}\sum_{i=0}^{p(k,m)-1} \mu \left (F^{-(i+\bar {p}(k,m))}
([u]_{-k})\cap R(k,m)\right)\eqno(1) 
$$

which is what we want. Now let us prove that the sequence $$ \left
(\frac{1}{n}\sum_{i=0}^{n-1}\mu (R(k,m))\cap F^{-i} ([u]_{-k}))\right
)_{m\in\NN} $$ converges uniformly in $n$ when $m\to\infty$.

We already know that for any real number $\varepsilon >0$,  for fixed $k$
there is an integer $m_0$ such that whenever $m\ge m_0$ one has $\mu
(R(k,m))\ge 1-\varepsilon$. Thus for any integer $i$ and
$m\ge m_0$ one has $$ \left \vert \mu ((X-R(k,m)) \cap F^{-i} ([u]_{-k}))\right
\vert \le \varepsilon , $$ hence $$ \left \vert \mu (F^{-i} ([u]_{-k}))- \mu
(R(k,m)\cap F^{-i} ([u]_{-k}))\right \vert \le \varepsilon . $$ For any integer
$n$ if $m\ge m_0$ one has $$ \left \vert \frac{1}{n}\sum_{i=0}^{n-1}\mu
(R(k,m)\cap F^{-i} ([u]_{-k}))- \frac{1}{n}\sum_{i=0}^{n-1}\mu ( F^{-i}
([u]_{-k}))\right\vert \le \frac {n\varepsilon}{n}= \varepsilon . $$

\medskip

Since the two convergence conditions hold, we have proved that the
two following limits exist and are the same:

$$\mu_c([u]_{-k})
=\lim_{n\to\infty}\lim_{m\to\infty}\frac{1}{n}\sum_{i=0}^{n-1}\mu \left
(R(k,m)\cap
F^{-i} ([u]_{-k})\right )
$$

$$
=\lim_{m\to\infty}\lim_{n\to\infty}\frac{1}{n} \sum_{i=0}^{n-1}\mu \left
(R(k,m)\cap
F^{-i} ([u]_{-k})\right ) =\lim_{n\to\infty}\mu_n ([u]_{-k}).
$$

Equality (1) permits to conclude that

$$
\mu_c([u]_{-k})=\lim_{m\to\infty} \frac{1}{p(k,m)}\sum_{i=0}^{p(k,m)-1} \mu
\left (R(k,m)\cap F^{-(i+\bar {p}(k,m))} ([u]_{-k})\right ). $$ \end{ppreuve}

The next corollary generalizes Theorem \ref{p7} to a larger class of \auts.
Its proof is straightforward.

\begin{cor} Let $(X,F,\mu )$ and $k \in \ZZ$ be such that $\mu $ is
$\sigma$--ergodic and equicontinuous for $(X,F\circ \sigma^{-k})$; then the
conclusions of Theorem \ref{p7} hold. \end{cor}


\subsection{The topological support of the measure $\mu_c$.} 

Remark first that the topological support $S(\mu_c)$ is contained in $W(F)$.

Recall that $R(k,m)$ is the set of points with at least one occurrence of $B$
between the coordinates $-k-m$ and $-m$ and another one between the coordinates
$m$ and $m+k$. We start with a technical lemma.

\begin{lem}\label{croi}
Let $\mu$ be a $\sigma$--ergodic measure, equicontinuous for $(X,F)$, and let
$B$ be a blocking word such that $\mu([B]_0)>0$. For any word $u$ in $A^{2k+1}$
the sequence $$ W_m(u)= \frac{1}{p(k,m)}\sum_{i=0}^{p(k,m)-1}
\mu (R(k,m)\cap F^{-(i+p'(k,m))} ([u]_{-k})) $$ is non--decreasing. 
\end{lem}

\begin{ppreuve}
Let $m_1<m_2$ be two natural integers. The two sequences $(U_i)=\left (\mu
(R(k,m_2)\cap F^{-i} ([u]_{-k}))\right )$ and $(V_i)=(\mu \left (R(k,m_1)\cap
F^{-i} ([u]_{-k}))\right )$ are ultimately periodic with preperiod and period
$p'(k,m_2)$ and $p(k,m_2)$ for the former, $p'(k,m_1)$ and $p(k,m_1)$ for
the latter.

Denote by $p'$ the greatest of the two integers $p'(k,m_2)$ and $p'(k,m_1)$,
and put $p=p(k,m_2) \times p(k,m_1)$. The sequences $(V_i)$ and
$(U_i)$ are ultimately periodic with preperiod $p'$ and period $p$
so that

$$
W_{m_2}(u)=\frac{1}{p}\sum_{i=0}^{p-1}
\mu \left (R(k,m_2)\cap F^{-(i+p')} ([u]_{-k}\right ). 
$$

Since $R(k,m_1)\subset R(k,m_2)$ one has
$$
W_{m_2}(u)\ge \frac{1}{p}\sum_{i=0}^{p-1} \mu \left (R(k,m_1)\cap F^{-(i+p')}
([u]_{-k}\right ) $$

$$
=\frac{1}{p(k,m_1)}\sum_{i=0}^{p(k,m_1)-1} \mu \left (R(k,m_1)\cap
F^{-(i+p'(k,m_1))} ([u]_{-k}\right ) =W_{m_1}(u). $$ 
\end{ppreuve}

\begin{pro}\label{Sincl}
Suppose $X$ is a transitive subshift of finite type, $F$ is onto and $\mu$
is equicontinuous for $(X,F)$ and $\sigma$--ergodic; then $S(\mu_c )\supset
S(\mu )$. 
\end{pro}

\begin{ppreuve}

Choose a blocking word  $B$ such that $\mu ([B]_0)>0$, and  $y \in
S(\mu)$; since $\mu([y(-k,k)]_{-k})>0$ for any integer  $k$ and
$\lim_{m\to\infty}\mu (R(k,m))=1$, there is an integer $m_0$ such that whenever
$m\ge m_0$ one has $\mu ([y(-k,k)]_{-k}\cap R(k,m))>0$.

For $\, m\in \NN$ choose a point $x_m$ in
$[y(-k,k)]_{-k}\cap R(k,m)$. By Proposition \ref{dense} the set of
$F$--periodic points is dense so there exists one, $y_m$, in the cylinder set
$[x_m(-k-m,k+m)]_{-k-m}$. The sequence $(F^n(x_m)(-k,k))$ does not depend on
the  coordinates to the left of $-k-m$ and to the right of $k+m$; it is
identical to the periodic sequence $F^n(y_m)(-k,k)$; in particular
$p'(k,m)=0$.

Fix $k$ and $m$: the sequence of sets
$(F^{-i} ([y(-k,k)]_{-k})\cap R(k,m))_{i\ge 0}$ is periodic; since $\mu
([y(-k,k)]_{-k}\cap R(k,m))>0$ one has, in the notation of Lemma \ref{croi},

$$
W_m(y(-k,k))=\frac{1}{p(k,m)}\sum_{i=0}^{p(k,m)-1}\mu (R(k,m)\cap F^{-i}
([y(-k,k)]_{-k})\ge 
$$

$$
\frac{1}{p(k,m)}\mu ((R(k,m))\cap [y(-k,k)]_{-k})>0. 
$$

By Proposition \ref{p7} and Lemma \ref{croi} the sequence $W_m(y(-k,k))$ is
non--decreasing and tends to $\mu_c ([y(-k,k)]_{-k})$ so that $\mu_c
([y(-k,k)]_{-k})>0$ and finally $y\in S(\mu_c )$.

\end{ppreuve}

In particular when $F$ is onto and $S(\mu )=A^\ZZ$ one has $S(\mu_c )=A^\ZZ$.


Let $(X,F)$ be a \autb having equicontinuity
points. For any blocking word $B$, call $E(F,B)$ the set of all points
$y\in X$ such that for any natural integer $k$, there is another natural integer
$m_0$ such that $\forall m\ge m_0$ and $\forall i\ge p'(k,m)$ one has
$[y(-k,k)]_{-k}\cap F^i(R(k,m))\neq \emptyset$.

\begin{pro}\label{mel}
The set $E(F,B)$ is a subshift; one has  $F(E(F,B)) \subset
E(F,B)$, thus $E(F,B)$ is contained in the limit set $W(F)$; if
$F$ is surjective, $E(F,B)=A^\ZZ$. If $X$ is transitive (resp.
strongly mixing) for $\sigma$, then  $E(F,B)$ does not depend on the choice of
the word  $B$ and can be denoted by $E(F)$; it is also transitive (resp.
strongly mixing) for $\sigma$. 
\end{pro}

\begin{ppreuve}
Since the definition of $E(F,B)$ depends only on the blocks of coordinates of
its points, $E(F)$ is a subshift; the fact that $F(E(F,B)) \subset
E(F,B)$ derives from the same remark. Let $X$ be transitive, $B$ and
$B'$ be two arbitrary blocking words; if $[y(-k,k)]_{-k}\cap F^i(R(k,m,B))\neq
\emptyset$, then $[y(-k,k)]_{-k}\cap F^i(R(k,m,B) \cap F^i(R(k,m',B'))\neq
\emptyset$ provided $m'$ is big enough, which implies
$E(F,B) \subset E(F,B')$. Transitivity or strong mixing result from the fact
that two words in $L(E(F))$ can occur in the image under $F^n$ of one
point containing one or several blocking words between their respective
occurrences. 
\end{ppreuve}

\begin{pro}\label{Eincl}
If $\mu$ is equicontinuous for $(X,F)$,  $S(\mu_c )\subset E(F)$. If moreover
$S(\mu )=X$, then $S(\mu_c )=E(F)$. \end{pro}

\begin{ppreuve}
Fix $B$ and assume that $y\not\in E(F)$, so there is an integer $k$ such that
for any integer $m$ and for any integer $i\ge p'(k,m)$ one has
$[y(-k,k)]_{-k}\not\subset F^i(R(k,m))$. Thus

$$
\frac{1}{p(k,m)}\sum_{i=0}^{p(k,m)-1}\mu (R(k,m)\cap
F^{-(i+p'(k,m))}([y(-k,k)]_{-k})) =0. 
$$
Applying Proposition \ref{p7} one obtains $\mu_c([y(-k,k)]_{-k})=0$ and
$y\not\in S(\mu_c )$. Hence $S(\mu_c )\subset E(F)$.

If $S(\mu )=X$ it is sufficient to prove that $S(\mu_c )\supset E(F)$.
If $x\in E(F)$, for any natural integer $k$ there is $m\in\NN$ such that the
union of cylinder sets $G=F^{-p'(k,m)}([x(-k,k)]_{-k})\cap R(k,m))$ is not
empty. Hence $\mu (G)>0$. By Lemma \ref{croi} and
Proposition \ref{p7} the sequence indexed by $m$

$$
\frac{1}{p(k,m)}\sum_{i=0}^{p(k,m)-1}\mu (R(k,m)\cap
F^{-(i+p'(k,m))}([y(-k,k)]_{-k}))
$$

is non--decreasing and since $\mu (G)>0$, Proposition \ref{p7} implies that $\mu_c$
 $([x(-k,k)]_{-k})$ $>0$ and finally $x\in S(\mu_c )$. \end{ppreuve}

{\bf Examples }\\
In \cite{Li86},  \cite{Ku94}, \cite{BM96} one can find examples of
\autsb that are surjective on $A^\ZZ$ and have equicontinuity
points. One can therefore apply
Proposition \ref{dense}, and  also Theorem
\ref{p7} and Proposition \ref{Sincl} if one assumes that  $\mu$ is for
instance a Bernoulli measure $B(p_1,p_2,p_3)$ different from the uniform
measure. In \cite{Mi} the automaton called ``Gliders and walls'' has
equicontinuity points without being onto.

Here is another example: the \aut $\, F: \{0,1,2\}^\ZZ \to \{0,1,2\}^\ZZ$ with
radius $1$ is defined by the local map $f$ such that $f(x_{-1},x_0,2)=x_0$,
$f(x_{-1},2,x_1)=2$ and when $x_1\in \{0,1\}$ and $x_0\neq 2$ then
$f(x_{-1},x_0,x_1)=x_0+x_1 \mbox{ mod } 2$. $F$ is onto; it has equicontinuity
points since $2$ is a blocking word. Let $\mu$ be
a Bernoulli measure with parameters $\{p,q,r\}$ on $A^\ZZ$; by
considering the cylinder sets $[2012]_0$ and $[2112]_0$ one easily checks that
the sequence $\mu \circ F^n$ does not converge vaguely. By
Theorem \ref{p7} it converges in Ces\`aro mean but, still considering the
same two cylinder sets, the limit cannot be the Bernoulli measure with
parameters $\{\frac{p+q}{2},\frac{p+q}{2},r\}$.

\medskip

{\bf Questions}

When is $\mu_c$ ergodic for $F$? When $F$ is onto and $\mu$ is the
uniform measure, which is $F$--invariant in this case, $\mu_c=\mu$ is never
$F$--ergodic (this would imply transitivity of $F$, which in its turn implies
sensitivity).

When is it $\sigma$--ergodic?

Are there conditions for  $\mu_c$ to be the uniform measure, or at least
Bernoulli or Markov? \bigskip

{\bf Acknowledgements}

We want to thank A. Maass for many suggestions; we are also grateful to the
referee, who signalled various shortcomings. Part of
this work was done by the second author at Universidad de Chile in Santiago,
thanks to the financial support of Fondap--Modelamiento Estocastico and
Ecos--Conicyt.

\end{document}